\documentclass[reqno]{amsart}
\usepackage{hyperref}
\usepackage{latexsym}

\newtheorem{theorem}{Theorem}[section]
\newtheorem{lemma}[theorem]{Lemma}

\newcommand{\bt}{\begin{theorem}}
\newcommand{\et}{\end{theorem}}
\newcommand{\blem}{\begin{lemma}}
\newcommand{\elem}{\end{lemma}}

\numberwithin{equation}{section}

\def\open#1{\setbox0=\hbox{$#1$}
\baselineskip = 0pt
\vbox{\hbox{\hspace*{0.4 \wd0}\tiny $\circ$}\hbox{$#1$}} 
\baselineskip = 11pt\!}

\def\fn{\open{f}}
\def\uno{\open{u}_1}
\def\unt{\open{u}_2}

\def\Pn{\open{P}}
\def\Rn{\open{R}}
\newcommand{\ep}{\hspace*{\fill}$\Box$}
\newcommand{\eps}{\varepsilon}

\newcommand{\R}{\mathbb R}
\newcommand{\N}{\mathbb N}

\newcommand{\om}{\omega}

\newcommand{\pa}{\partial}

\newcommand{\CC}{\mathcal{C}}

\newcommand{\nn}{\nonumber}
\newcommand{\beq}{ \begin{equation} }
\newcommand{\eeq}{\end{equation} }
\newcommand{\bea}{\begin{eqnarray}}
\newcommand{\eea}{\end{eqnarray}}
\newcommand{\beas}{\begin{eqnarray*}}
\newcommand{\eeas}{\end{eqnarray*}}
\newcommand{\beqs}{\begin{equation*}}
\newcommand{\eeqs}{\end{equation*}}

\newcommand{\ben}{\begin{enumerate}}
\newcommand{\een}{\end{enumerate}}
\newcommand{\ba}{\begin{array}}
\newcommand{\ea}{\end{array}}
\newcommand{\brem}{\begin{thr} {\bf Remark. }\rm}
\newcommand{\ethi}{\end{thr}}

\begin{document}

\title[Classical Solutions of the VKG System]%
{On Classical Solutions of the Relativistic Vlasov-Klein-Gordon System}

\author[M. Kunzinger]{Michael Kunzinger}
\address{Fakult\"at f\"ur Mathematik, Universit\"at  Wien\\
Nordbergstrasse 15\\ 1090 Wien\\ Austria}
\email{Michael.Kunzinger@univie.ac.at}
\urladdr{http://www.mat.univie.ac.at/\~{}mike/}
\author[G. Rein]{Gerhard Rein}
\address{Fakult\"at f\"ur Mathematik und Physik\\
Universit\"at Bayreuth\\
95440 Bayreuth\\ Germany}
\email{gerhard.rein@uni-bayreuth.de}
\urladdr{http://www.uni-bayreuth.de/departments/math/org/mathe6/staff/memb/grein/home.html}
\author[R. Steinbauer]{Roland Steinbauer}
\address{Fakult\"at f\"ur Mathematik, Universit\"at Wien\\
Nordbergstrasse 15\\ 1090 Wien\\ Austria}
\email{Roland.Steinbauer@univie.ac.at}
\urladdr{http://www.mat.univie.ac.at/\~{}stein/}
\author[G. Teschl]{Gerald Teschl}
\address{Fakult\"at f\"ur Mathematik,  Universit\"at Wien\\
Nordbergstrasse 15\\ 1090 Wien\\ Austria\\ and International Erwin Schr\"odinger
Institute for Mathematical Physics, Boltzmanngasse 9\\ 1090 Wien\\ Austria}
\email{Gerald.Teschl@univie.ac.at}
\urladdr{http://www.mat.univie.ac.at/\~{}gerald/}
\thanks{Partially supported by the Austrian Science Fund's 
Wittgenstein 2000 Award of P.~A.~Markowich.}

\begin{abstract}
We consider a collisionless ensemble of classical particles coupled to a
Klein-Gordon field. For the resulting nonlinear system
of partial differential equations, the relativistic
Vlasov-Klein-Gordon system, we prove local-in-time existence of classical
solutions and a continuation criterion which says that a solution can blow up
only if the particle momenta become large. We also show that classical solutions are
global in time in the one-dimensional case.
\end{abstract}

\maketitle

\section{Introduction}\label{intro}
\setcounter{equation}{0}

In kinetic theory one often considers collisionless ensembles of 
classical particles which interact only by fields which they create
collectively. This situation is commonly referred to as the mean field
limit of a many particle system. Such systems have been  studied extensively. 
In the case of non-relativistic, gravitational
or electrostatic fields the corresponding system of partial differential
equations is the Vlasov-Poisson system, in the case of relativistic electrodynamics
it is the Vlasov-Maxwell system and in the case of general relativistic gravity
the Vlasov-Einstein system.

On the other hand the coupling of a single particle to a
classical or quantum field has been studied. In case of the Maxwell field this is a
classical problem \cite{Ab}, but the actual dynamics and asymptotics
of such systems is still an active area of research \cite{IKS,KKS1,KKS2,KS}.
In \cite{IKM} the case of a single classical particle coupled to a
quantum mechanical Klein-Gordon field was investigated. 

In the present paper we consider a collisionless ensemble of particles 
moving at relativistic speeds, coupled to a Klein-Gordon field. This is
a natural generalization of the one-particle situation just described.
Let $f=f(t,x,v)\geq 0$ denote the density of the particles
in phase space, $\rho=\rho(t,x)$ their density in space, and
$u=u(t,x)$ a scalar Klein-Gordon field; $t \in \R$, $x \in \R^3$, and
$v \in \R^3$ denote time, position, and momentum, respectively.
The system then reads as follows: 
\begin{equation} \label{vl}
\pa_t f + \hat v \cdot \pa_x f - \pa_x u \cdot \pa_v f = 0, 
\end{equation}
\begin{equation} \label{kg}
\pa_t^2 u - \Delta u + u = - \rho, 
\end{equation}
\begin{equation} \label{rhodef}
\rho(t,x) = \int f(t,x,v)\,dv.
\end{equation}
Here we have set all physical constants as well as the rest mass
of the particles to unity, and
\begin{equation} \label{hatvdef}
\hat v = \frac{v}{\sqrt{1+|v|^2}}
\end{equation}
denotes the relativistic velocity of a particle with momentum $v$.
This system is called the relativistic Vlasov-Klein-Gordon system.
It is supplemented by 
initial data
\begin{equation} \label{id}
f(0)= \fn,\ u(0) = \uno,\ \pa_t u(0) = \unt .
\end{equation}

The study of this system was initiated in \cite{KRST} where the existence of
global weak solutions for initial data satisfying a size restriction was proved.
This size restriction was necessary because  
the energy of the system is indefinite so that conservation of energy does not lead 
to a-priory bounds for general data. 
A natural next step in the study of the Vlasov-Klein-Gordon
system is the existence theory of classical solutions, locally and if possible
globally in time. This is the topic of the present investigation.

Another motivation for studying this system of partial differential equations 
is an intrinsically mathematical one. Since the field equation is hyperbolic 
the system resembles
the relativistic Vlasov-Maxwell system, for which the quest for global-in-time
classical general solutions is still open. One might hope that studying
related systems can help in understanding these open problems more thoroughly. 
In fact in this work we follow the general outline of the existence proof of 
Glassey and Strauss \cite{GlStr2} for the Vlasov-Maxwell system. 
Note, however, that the existence theory of weak solutions of the two systems 
is quite different \cite{DL,KRST,R}. 

The paper proceeds as follows. In Section~\ref{estimates} we prove some
a-priori estimates necessary for the proof of our main result. These estimates
rely on representation formulas for the first and second order derivatives of 
the Klein-Gordon field $u$, cf.\ Lemmas~\ref{lem1derivative} and \ref{lem2derivative}. 
Note that in contrast to the
corresponding parts in \cite{GlStr2} we also need to bound the mixed second
order derivatives of the field. In Section~\ref{secex} we prove
our main results, a local-in-time existence and uniqueness result for classical 
solutions and a continuation criterion
which says that such solutions can blow up in finite time only if the support of $f$
in momentum space becomes unbounded,
cf.\ Thms~\ref{locsol} and ~\ref{contcrit}. 
In Section~\ref{1d} we briefly show that the continuation criterion
is indeed satisfied in the one-dimensional situation where $x,v \in \R$
so that we obtain global classical solutions in that case.
Finally some material on the (inhomogeneous) 
Vlasov equation is collected in an appendix as well as some unpleasant
technical aspects of the proof of the local existence result, which often have been 
omitted in the treatment of related systems.
 
\section{A-priori estimates}\label{estimates}
\setcounter{equation}{0}

Although our notation is mostly standard or self-explaining
we mention the following conventions:
For a function $h=h(t,x,v)$ or $h=h(t,x)$ we denote for given $t$
by $h(t)$ the corresponding function of the remaining variables.
For a function $h$ depending on the variables $x,v$ we denote its gradient
by $\pa_{(x,v)}h$. By $\|\,.\,\|_p$ we denote the usual $L^p$-norm for 
$p\in[1,\infty]$. 
The subscript $c$ in function spaces refers to compactly supported
functions. Sometimes we write $z = (x,v) \in \R^3 \times \R^3$.

One main ingredient of our analysis are
representation formulas for $u$ and its derivatives, which will allow us to
establish the necessary a-priori bounds.
As point of departure we recall that the solution of (\ref{kg}) is given by
\begin{equation} \label{solkg}
u(t,x) = u_{\mathrm{hom}}(t,x) + u_{\mathrm{inh}}(t,x),\quad t\geq 0,\ x\in \R^3,
\end{equation}
where
\begin{eqnarray*}
u_{\mathrm{hom}}(t,x)
&=&
\frac{1}{4 \pi t^2} \int_{|x-y|=t} \uno(y)\, dS_y
- \frac{1}{4 \pi t^2} \int_{|x-y|=t} (\pa_x \uno)(y)\cdot y\, dS_y\\ 
&&
{}
- \frac{1}{8 \pi} \int_{|x-y|=t} \uno(y)\, dS_y
- \frac{1}{4 \pi} \int_{|x-y|\leq t} \uno(y)\, 
\left(\frac{J_1(\xi)}{\xi}\right)'
\frac{t}{\xi} dy\\ 
&&
{}
+ \frac{1}{4 \pi t} \int_{|x-y|=t} \unt(y)\, dS_y
- \frac{1}{4 \pi} \int_{|x-y|\leq t} \unt(y)\, 
\frac{J_1(\xi)}{\xi} dy,
\end{eqnarray*}
is the solution of the homogeneous Klein-Gordon equation with initial data as 
in (\ref{id}) and $\xi:= \sqrt{t^2 - |x-y|^2}$, and
\begin{eqnarray*}
u_{\mathrm{inh}}(t,x)
&=& 
- \frac{1}{4 \pi}\int_0^t \int_{|x-y|=t-s}\!\!\!\!\!\!
\rho(s,y) \,dS_y\frac{ds}{t-s} 
+ \frac{1}{4 \pi}\int_0^t \int_{|x-y|\le t-s}\!\!\!\!\!\! 
\rho(s,y)\,\frac{J_1(\xi)}{\xi} dy\,ds
\end{eqnarray*}
with $\xi:=\sqrt{(t-s)^2 - |x-y|^2}$
is the solution of the inhomogeneous Klein-Gordon equation with vanishing
initial data, cf.\ \cite{MS} or \cite{Sid}; $J_k$ denotes the Bessel function.
To derive formulas for the derivatives of $u$
the differential operators
\[
S = \pa_t + \hat{v} \cdot \pa_x,\
T = - \om \pa_t + \pa_x; \qquad \om=\frac{x-y}{|x-y|},
\]
which are adapted to our system and have first been introduced in \cite{GlStr2}
in connection with the Vlasov-Maxwell system turn out to be useful.

\begin{lemma}\label{lem1derivative} (Representation of $\pa u$)\\
Suppose $u\in \CC^2$ is a solution of the Klein-Gordon equation (\ref{kg}) with
$\rho$ given by (\ref{rhodef}) for some $f\in\CC^1$. Then
\[
\pa_k u(t,x) = F_0^k(t,x) +  F_S^k(t,x) + F_T^k(t,x) + F_R^k(t,x) + F_J^k(t,x),
\quad k\in \{1,2,3,t\}
\]
where $F_0^k$ is a linear functional of the initial data only, and
\beas
F^k_S(t,x) &=& -\frac{1}{4\pi} \int_{|x-y|\leq t}\int\frac{\om_k}{1+\om\cdot\hat{v}}
(S f)(t-|x-y|,y,v) \,dv \frac{dy}{|x-y|},\\
F^k_T(t,x) &=& \frac{1}{4\pi} \int_{|x-y|\leq t}\int a^k(\om,\hat{v})
f(t-|x-y|,y,v) \,dv \frac{dy}{|x-y|^2},\\
F^k_R(t,x) &=& \frac{1}{8\pi} \int_{|x-y|\le t}
\rho(t-|x-y|,y) \om_k dy,\\
F^k_J(t,x) &=& -\frac{1}{4\pi} \int_0^t \int_{|x-y|\le t-s} \rho(s,y) 
\frac{J_2(\xi)}{\xi^2}
(x_k-y_k) \,dy\,ds,\\
F^t_S(t,x) &=& -\frac{1}{4\pi} \int_{|x-y|\leq t}\int \frac{1}{1+ \om\cdot\hat{v}}
(S f)(t-|x-y|,y,v) \,dv \frac{dy}{|x-y|},\\
F^t_T(t,x) &=& \frac{1}{4\pi} \int_{|x-y|\leq t}\int a^t(\om,\hat{v})
f(t-|x-y|,y,v) \,dv \frac{dy}{|x-y|^2},\\
F^t_R(t,x) &=& \frac{1}{8\sqrt{2}\pi} \int_{|x-y|\le t} \rho(t-|x-y|,y)
dy,\\
F^t_J(t,x) &=& -\frac{1}{4\pi} \int_0^t \int_{|x-y|\le t-s} \rho(s,y) 
\frac{J_2(\xi)}{\xi^2} (t-s)
dy ds.
\eeas
Here $\xi=\sqrt{(t-s)^2-|x-y|^2}$ and the kernels $a^k$ and $a^t$ are
\[
a^k(\om,\hat{v}) = - \frac{\hat{v}_k}{1+ \om\cdot\hat{v}} - 
\frac{\om_k}{(1+|v|^2)(1+ \om\cdot\hat{v})^2}, \qquad
a^t(\om,\hat{v}) = \frac{|\hat{v}|^2 + \om \cdot \hat{v}}{(1+ \om\cdot\hat{v})^2}.
\]
\end{lemma}

\begin{proof}
The proof is a straightforward calculation using
\[
\pa_{x_k} = \frac{\om_k}{1+ \om \cdot \hat{v}}\, S + 
\sum_{j=1}^3 \left( \delta_{jk} -
\frac{\om_k \hat{v}_j}{1+ \om \cdot \hat{v}} \right)\, T_j,\quad
\pa_t  = \frac{1}{1+ \om \cdot \hat{v}} \, ( S - \hat{v} \cdot T )
\]
where $\delta_{jj}=1$ and $\delta_{jk}=0$ for $j\ne k$. Terms involving only the
initial data are collected in $F_0^k$ and terms involving
$T f$ are integrated by parts, using the identity
$(T f)(t-|x-y|,y,v) = \pa_y (f(t-|x-y|,y,v))$.
\end{proof}

Next we turn to the second order derivatives.

\begin{lemma}\label{lem2derivative} (Representation of $\pa^2 u$)\\
Suppose $u\in \CC^2$ is a solution of the Klein-Gordon equation (\ref{kg}) with
$\rho$ given by (\ref{rhodef}) for some $f\in\CC^2$. Then we have 
for $k,\ell\in\{1,2,3,t\}$
\beas
\pa_{k\ell} u(t,x) &=& F_0^{k\ell} + F_{SS}^{k\ell} + F_{ST}^{k\ell} + 
F_{TS}^{k\ell} + F_{TT}^{k\ell} + F_{RS}^{k\ell} + F_{RT}^{k\ell} + F_{JR}^{k\ell} + 
F_{JJ}^{k\ell},
\eeas
where $F_0^{k\ell}$ are linear functionals of the initial data only, and
\beas
F_{SS}^{k\ell} &=& \frac{-1}{4\pi} \int_{|x-y|\leq t}\int c^{k\ell}(\om,\hat{v})
(S^2 f)(t-|x-y|,y,v) \,dv \frac{dy}{|x-y|},\quad
|c^{k\ell}| \le \frac{C}{(1+\om\cdot\hat{v})^2}\\
F_{ST}^{k\ell} &=& \frac{1}{4\pi} \int_{|x-y|\leq t}\int b_1^{k\ell}(\om,\hat{v})
(S f)(t-|x-y|,y,v) \,dv \frac{dy}{|x-y|^2},\quad
|b_1^{k\ell}| \le \frac{C}{(1+\om\cdot\hat{v})^3}\\
F_{TS}^{k\ell} &=& \frac{1}{4\pi}  \int_{|x-y|\leq t}\int b_2^{k\ell}(\om,\hat{v})
(S f)(t-|x-y|,y,v) \,dv \frac{dy}{|x-y|^2}, \quad
|b_2^{k\ell}| \le \frac{C}{(1+\om\cdot\hat{v})^3}\\
F_{TT}^{k\ell} &=& \frac{-1}{4\pi} \int_{|x-y|\leq t}\int a^{k\ell}(\om,\hat{v})
f(t-|x-y|,y,v) \,dv \frac{dy}{|x-y|^3}, \quad
\int_{|\om|=1} a^{k\ell}(\om,\hat{v}) d\om =0,\\
F_{RS}^{k\ell} &=& \frac{1}{8\pi} \int_{|x-y|\leq t}\int d^{k\ell}(\om,\hat{v})
(S f)(t-|x-y|,y,v) \,dv \,dy,\quad
|d^{k\ell}| \le \frac{C}{1+\om\cdot\hat{v}}\\
F_{RT}^{k\ell} &=& \frac{1}{8\pi}  \int_{|x-y|\leq t}\int e^{k\ell}(\om,\hat{v})
f(t-|x-y|,y,v) \,dv \frac{dy}{|x-y|},\quad
|e^{k\ell}| \le \frac{C}{(1+\om\cdot\hat{v})^2}\\
F_{JR}^{k\ell} &=&  \frac{-1}{32\pi}
\big(1-\frac{\sqrt{2}-1}{\sqrt{2}}
(\delta_{kt}+\delta_{\ell t}-\delta_{kt}\delta_{\ell t})\big) 
\int_{|x-y|\le t} \rho(t-|x-y|,y) 
\om_k \om_\ell |x-y| \,dy,\\
F_{JJ}^{k\ell} 
&=&
\frac{2(\delta_{kt}+\delta_{\ell t}-\delta_{kt}\delta_{\ell t})-1}{4\pi}\\
&& \int_0^t 
\int_{|x-y|\le t-s} \!\!\!\! \rho(s,y) \Big(
\frac{J_3(\xi)}{\xi^3}(x_k-y_k)(x_\ell-y_\ell) + (1-\delta_{kt}-\delta_{\ell t})
\frac{J_2(\xi)}{\xi^2}\Big)
dy\,ds.
\eeas
Here we have set $x_t=t$, $y_t=s$, and $\om_t=1$ for notational convenience.
\end{lemma}

\begin{proof}
The proof is a long and tedious calculation similar to the previous lemma. 
The critical part is, of course, to prove that the kernels $a^{k\ell}$ appearing
in the singular $TT$-terms vanish when integrated over the unit sphere. To give some
flavor of the respective calculations we outline them in the case of the 
(most complicated) kernel $a^{k\ell}$ for $1\leq k,\ell\leq 3$. We have
\begin{eqnarray*}
 a^{k\ell}&=&-3\,\frac{\om_\ell\cdot\hat v_k+\om_k\cdot\hat v_\ell}
                  {(1+|v|^2)(1+\om\cdot\hat v)^3}
            -\frac{3\om_k\om_\ell}
                  {(1+|v|^2)^2(1+\om\cdot\hat v)^4}\\
            &&{}-\frac{2\hat v_k\hat v_\ell}
                  {(1+\om\cdot\hat v)^2}
            +\frac{\delta_{k\ell}}
                  {(1+|v|^2)(1+\om\cdot\hat v)^2}\\
          &=:&a^{k\ell}_1+a^{k\ell}_2+a^{k\ell}_3+a^{k\ell}_4.      
\end{eqnarray*}
Using $\pa_{v_i}\big((\sqrt{1+|v|^2}+\om\cdot v)^{-2}\big)=-2(\hat v_i+\om_i)\,
(\sqrt{1+|v|^2}+\om\cdot v)^{-3}$ we find
\begin{eqnarray*}
  v_k\int_{|\om|=1}\frac{\om_\ell\,d\om}{(\sqrt{1+|v|^2}+\om\cdot v)^3}
  &=&-\frac{v_k}{2}\pa_{v_\ell}
      \int_{|\om|=1}\frac{d\om}{(\sqrt{1+|v|^2}+\om\cdot v)^2}\\
  &&-v_k\hat v_\ell\int_{|\om|=1}
      \frac{d\om}{(\sqrt{1+|v|^2}(1+\om\cdot\hat v))^3}
  \,=\,-4\pi v_k v_\ell.
\end{eqnarray*}
Hence $\int_{|\om|=1}a^{k\ell}_1\,d\om=24\pi v_kv_\ell$.
By the identity
\[
 \pa_{v_k}\Big(\pa_{v_\ell}\frac{1}{(\sqrt{1+|v|^2}+\om\cdot v)^2}
               +\frac{2\hat v_\ell}{(\sqrt{1+|v|^2}+\om\cdot v)^3}
         \Big)
 = 6\,\frac{\om_\ell(\hat v_k+\om_k)}{(\sqrt{1+|v|^2}+\om\cdot v)^4}
\] 
we obtain
\begin{eqnarray*}
  -\int_{|\om|=1}\frac{\om_\ell\hat v_k\,d\om}{(\sqrt{1+|v|^2}+\om\cdot v)^4}\\
 &&\hspace*{-1cm}
  +\frac{1}{6}\pa_{v_k}
     \int_{|\om|=1}\Big(\pa_{v_\ell}\frac{1}{(\sqrt{1+|v|^2}+\om\cdot v)^2}
                                +\frac{2\hat v_\ell}{(\sqrt{1+|v|^2}+\om\cdot v)^3}
                          \Big)\,d\om\\
 &&\hspace*{-1.2cm}
 =\frac{16\pi}{3}v_\ell v_k+\frac{4\pi}{3}\delta_{\ell k},
\end{eqnarray*}
which implies 
$\int_{|\om|=1}a^{k\ell}_2\,d\om=-16\pi v_kv_\ell-4\pi\delta_{k\ell}$.
The remaining two terms can be integrated directly to yield
\[
 \int_{|\om|=1}a^{k\ell}_3\,d\om =-8\pi v_k v_\ell,\qquad
 \int_{|\om|=1}a^{k\ell}_4\,d\om =4\pi\delta_{k\ell}.
\]
Summing up we obtain the desired result
\[
 \int_{|\om|=1}a^{k\ell}\,d\om=0.
\]
The computations for the kernels 
\begin{eqnarray*}
 a^{tk}&=&\frac{2\hat v_k(\om\cdot\hat v+|\hat v|^2)}{(1+\om\cdot\hat v)^3}
         +\frac{3\om_k(\om\cdot\hat v+|\hat v|^2)}{(1+|v|^2)(1+\om\cdot\hat v)^4}
         -\frac{\hat v_k}{(1+|v|^2)(1+\om\cdot\hat v)^3}\qquad( 1\leq k\leq 3) \\
 a^{tt}&=&\frac{1}{(1+\om\cdot\hat v)^4}
          \big(3|\hat v|^4-(\om\cdot\hat v)^2|\hat v|^2-|\hat v|^2+3(\om\cdot\hat v)^2
               +4\om\cdot\hat v|\hat v|^2
          \big)
\end{eqnarray*}
proceed along the same lines.
\end{proof}

Note that in contrast to the corresponding analysis of the Vlasov-Maxwell
system \cite{GlStr2} we also need to provide the mixed second order 
derivatives of the field. The fact that the
averaging property of corresponding the $TT$-kernel holds is a pleasant surprise 
in so far as the mixed derivatives do not appear in the field equation.

We can now use these formulas to derive the necessary estimates on the derivatives
of $u$. 
Constants denoted by $C$ are positive and may change their value from line to line.
If they depend on the initial data this is explicitly mentioned.

\begin{lemma} \label{lemfuestimates} 
(Estimates on $\pa f$, $\pa u$, $\pa^2 u$)
\begin{itemize}
\item[(i)] 
Suppose that $f\in\CC^1$ is a solution of the Vlasov equation
\[
S f(t,x,v) = F(t,x) \pa_v f(t,x,v), \qquad f(0,x,v) = \fn(x,v), \qquad t\in [0,T[,
\]
for some $F\in\CC^1$. Then we have, with $z=(x,v)$,
\[
\| \pa_z f(t) \|_\infty \le \| \pa_z \fn \|_\infty +
C \int_0^t (1 + \| \pa _x F(s) \|_\infty) \| \pa_z f(s) \|_\infty ds.
\]
\item[(ii)] 
In addition, assume
there exists an increasing function $P(t)$ such that $f(t,x,v)=0$ for 
$|v|\ge P(t)$, and suppose
$u$ is the solution (\ref{solkg}) of the Klein-Gordon equation
for $t\in [0,T[$. Then $u\in\CC^2$, and
\[
\| \pa_{(t,x)} u(t) \|_\infty \le C 
\left((1+t)^5 (1+P(t))^5  +  \int_0^t  (t-s)(1+P(s))^6 \| F(s) \|_\infty ds \right),
\]
where the constant $C$ depends on the norms of the initial data. 
\item[(iii)]
On any bounded time interval on which $F$ and $P$ are bounded we have
\[ \qquad\qquad
 \| \pa_{(t,x)}^2 u(t) \|_\infty 
\le 
C \left( 1 +
 \log^* (\sup_{0\le s \le t} \| \pa_x f(s) \|_\infty )
 +  \int_0^t  \| \pa_{(t,x)} F(s) \|_\infty ds \right),
\]
where $\log^*(s)=s$ for $0\leq s \leq 1$ and
$\log^*(s)=1 + \log(s)$ for $1 \ge s$, and $C$ depends on 
the time interval, the bounds for $P$ and $F$, and the initial data.
\end{itemize}
\end{lemma}

\begin{proof} The assertion in (i) follows directly from (\ref{inteqfz})  below
with $g = 0$ and $G(t,z)=(\hat v,F(t,x))$.
As to (ii), the estimate of $\pa u$ is a straightforward consequence 
of Lemma~\ref{lem1derivative}
by using $S f = F \pa_v f$ and getting rid of the $v$ derivatives via integration
by parts. Then we estimate each term individually using 
$\| f(t) \|_\infty = \| \fn \;\|_\infty$, cf.\ (\ref{solvl}) with $g\equiv 0$, 
the estimate
\[
 \frac{1}{1+ \om \cdot \hat{v}} \le 2(1+|v|^2),
\]
and the fact that the $v$-integration is only over a finite ball
of radius $P(t)$.

To prove (iii) we use the same procedure with Lemma~\ref{lem2derivative}
replacing Lemma~\ref{lem1derivative}. First we need to get rid of the
second derivative (i.e., the $S^2 f$ term in $F_{SS}^{k\ell}$). So let us first assume
$S f \in \CC^1$  and consider
\beas
F_{SS}^{k\ell} &=& \frac{-1}{4\pi} \int_{|x-y|\leq t}\int \int c^{k\ell}(\om,\hat{v})
(S^2 f)(t-|x-y|,y,v) \,dv \frac{dy}{|x-y|}.
\eeas
Observe that
\[
S^2 f = (S F)\cdot \pa_v f + F \cdot \pa_v (F \cdot \pa_v f) - F c(v) \pa_x f,
\quad c_{jk}(v)= \pa_{v_j} \hat{v}_k.
\]
Inserting this into $F_{SS}^{k\ell}$ and removing all $v$ derivatives using
integration by parts we end up with an expression involving only first order derivatives
of $f$, which also holds without the $S f \in \CC^1$ assumption by a standard
approximation argument, in particular, $u\in\CC^2$.

Now we can estimate each term as before except for the $TT$-kernel which is the
most critical one, i.e.,
\beas
 F_{TT}^{k\ell} 
 &=&\frac{1}{4\pi} \int_{|x-y|\leq t}\int a^{k\ell}(\om,\hat{v})
     f(t-|x-y|,y,v) \,dv \frac{dy}{|x-y|^3}\\
 &=&\int_0^t \frac{1}{t-s} \int_{|\om|=1} 
     \int a^{kl}(\om, \hat{v}) f(s,x+(t-s)\om,v)dv\,d\om\,ds,
\eeas
where we split the $s$-integral into two parts, $I$ over $[0,t-\tau]$ and $II$ over 
$[t-\tau,t]$. The first of these can be estimated directly by
$C\log (t/\tau)$.
For the second one we use the averaging property of $a^{k\ell}$ to write it as
\[
 II=\int_{t-\tau}^t\frac{1}{t-\tau}\int_{|\om|=1}\int
 a^{k\ell}(\om,\hat v)\Big(f(s,x+(t-s)\om,v)-f(s,x,v)dvd\om\Big) ds.
\]
Hence by the mean value theorem
\[
 |II|\leq C\ \tau\sup_{t-\tau\leq s\leq t}\|\pa_x f(s)\|_\infty.
\]
Summing up we obtain the estimate
\[
| F_{TT}^{k\ell}(t) | \le C ( \log(t/\tau) + \tau N(t) ), \quad 0\le t \le T,
\]
where $N(t) = \sup_{0\le s \le t} \| \pa_x f(s) \|_\infty$.
For $N(t) \le t^{-1}$ we can choose the optimal value $\tau= N(t)^{-1}$, otherwise
we choose $\tau= t$, and this yields
\[
| F_{TT}^{k\ell}(t) | \le C \log^*(t N(t)).
\]
Combining these estimates the remaining claim follows.
\end{proof}

\section{Existence of classical solutions}\label{secex}
\setcounter{equation}{0}

We now have collected all ingredients to show existence and uniqueness of local
classical solutions of the relativistic Vlasov-Klein-Gordon system.

\begin{theorem}\label{locsol}
(Local existence of classical solutions)\\
Let $\fn\; \in \CC_c^1(\R^6)$, $\uno \in \CC_b^3(\R^3)$,
$\unt \in \CC_b^2(\R^3)$.
Then there exists a unique classical solution 
\[ 
f \in \CC^1([0,T[ \times \R^6)),\ u\in \CC^2([0,T[\times\R^3)
\]
of the relativistic Vlasov-Klein-Gordon system (\ref{vl})--(\ref{rhodef})
for some $T>0$, satisfying the initial conditions (\ref{id}).
Moreover, 
\[
f(t,x,v) = 0\ \mbox{for}\ |x| \geq \Rn + t \ \mbox{or}\ |v| \geq P(t)
\]
where $\open{R}$ is determined by $\fn$ and $P$ is a positive continuous function
on $[0,T[$.
\end{theorem}

\begin{proof}
We begin with the {\em uniqueness} part which relies 
only on Lemma~\ref{lem1derivative}. Let $(f^{(1)},u^{(1)})$, 
$(f^{(2)},u^{(2)})$ be two solutions satisfying the same
initial conditions, and on any compact time interval
$[0,T_0]$ on which both solutions exist define
$f=f^{(1)} -f^{(2)}$ and  $u=u^{(1)} -u^{(2)}$.
Then $u$ satisfies (\ref{kg}) with (\ref{rhodef}), and $f$ satisfies
\[
S f = \pa_x u^{(1)} \pa_v f^{(1)} - \pa_x u^{(2)} \pa_v f^{(2)} =
\pa_x u \pa_v f^{(1)} + \pa_x u^{(2)} \pa_v f.
\]
Proceeding on a term-by-term basis using the representation from 
Lemma~\ref{lem1derivative} and replacing the $Sf$ term via the Vlasov equation
and integrating by parts we obtain the estimate
\[
\| \pa_x u(t) \|_\infty 
\le C \int_0^t \left(\|f(s)\|_\infty +
\| \pa_x u(s)\|_\infty \| f^{(1)}(s) \|_\infty + 
\| \pa_x u^{(2)}(s)\|_\infty \| f(s) \|_\infty \right)\, ds.
\]
Using the boundedness of $f^{(1)}$ and $\pa_x u^{(2)}$  this implies
\[
 \| \pa_x u(t) \|_\infty \le C
 \int_0^t ( \|f(s)\|_\infty +  \| \pa_x u(s) \|_\infty )\, ds.
\]
On the other hand, by (\ref{solvl})
\[
f(s,x,v)=  \int_0^t (\pa_x u \,\pa_v f^{(1)}) (s,Z^{(2)}(s,t,x,v))\, ds,
\]
where $Z^{(2)}(s,t,x,v)$ is the solution of the characteristic
system corresponding to $u^{(2)}$. Hence we obtain, using the boundedness 
of $\pa_v f^{(1)}$,
\[
\| f(t) \|_\infty \le C \int_0^t \| \pa_x u(s) \|_\infty ds.
\]
Combining both estimates gives
\[
\| f(t) \|_\infty + \| \pa_x u(t) \|_\infty \le C 
\int_0^t ( \|f(s)\|_\infty + \| \pa_x u(s) \|_\infty )\, ds
\]
and Gronwall's lemma implies $f(t)=u(t)=0$, proving uniqueness.
\medskip

Next we turn to {\em existence}. To this end we set up an iterative scheme and
prove its convergence to a solution. Let $f^{(0)}=\fn$, $u^{(0)}= \uno$ and
define $f^{(n)}\in \CC^1([0,\infty[ \times\R^6))$,
$u^{(n)}\in\CC^2([0,\infty[ \times\R^3)$ 
recursively via
\[
S f^{(n)} - \pa_x u^{(n-1)} \pa_v f^{(n)} =0, \quad f^{(n)}(0)=\fn
\]
and
\[
\pa_t^2 u^{(n)} - \Delta u^{(n)} + u^{(n)} = - \rho^{(n)}, 
\quad u^{(n)}(0)=\uno, \quad \pa_tu^{(n)}(0)=\unt
\]
where
\[
\rho^{(n)}(t,x) = \int f^{(n)}(t,x,v)\,dv;
\]
notice that $f^{(n)}$ satisfies a support estimate as the one 
asserted for the solution
$f$, but for all $t\geq 0$ and with the increasing function 
\[
P^{(n)}(t):=\sup\left\{|v| \mid
f^{(n)}(\tau,x,v)\not=0,\ 0\leq\tau\leq t,\ x\in\R^3\right\}
\]
instead of $P$. 

\noindent
{\em Step 1 } (Uniform bounds on $P^{(n)}$ and $\pa u^{(n)}$):\\
By Lemma~\ref{lemfuestimates} (ii)  and (\ref{solvl}) we have
\[
\| \pa_{(t,x)} u^{(n)}(t) \|_\infty 
\le 
C \left( (1+t)^5 (1+P^{(n)})^5 + 
\int_0^t s^2 P^{(n)}(s)^6 \| \pa_{(t,x)}u^{(n-1)}(s) \|_\infty ds\right)
\]
and
\[
P^{(n)}(t) \le \Pn + \int_0^t  \| \pa_x u^{(n-1)}(s) \|_\infty ds,
\]
where $\Pn$ bounds the $v$-support of the data.\\
Define $Q^{(n)}(t):= \max_{0\leq k\leq n} \| \pa_{(t,x)} u^{(k)}(t) \|_\infty$.
Then
\[
Q^{(n)}(t)
\leq
C \left( (1+t)^5 (1+P^{(n)}(t))^5 + 
\int_0^t (t-s)(1+ P^{(n)}(s))^6  Q^{(n)}(s)\, ds\right),
\]
and by Gronwall's lemma,
\[
Q^{(n)}(t)
\leq
C \left( (1+t)^5 (1+P^{(n)}(t))^5 
\exp\big(t^2(1+P^{(n)}(t))^6 \big)\right).
\]
If we insert this into the estimate for $P^{(n)}$ we find that
\beq \label{pngronwall}
P^{(n)}(t) \le \Pn + 
C \int_0^t (1+s)^5 (1+P^{(n)}(s))^5 
\exp\left(s^2(1+P^{(n)}(s))^6\right)ds,\ t\geq 0.
\eeq
Hence by induction, $P^{(n)}(t) \leq P(t),\ n\in \N,\ t\in [0,T[$,
where $P$ is the maximal solution of the integral equation
corresponding to (\ref{pngronwall}), which exists on some time
interval $[0,T[$ whose length is determined by $\Pn$ and the norms
of the initial data entering the constant $C$.
In addition, $\pa_{(t,x)}u^{(n)}(t)$ are bounded in terms of $P(t)$
on that interval.

For the rest of the proof we now argue on a bounded, arbitrary, but 
fixed time interval 
$[0,T_0] \subset [0,T[$. Constants denoted by $C$ may now depend on $T_0$
and the bounds established in Step~1.

\noindent
{\em Step 2 } (Uniform bounds on $\pa^2 u^{(n)}$ and $\pa_z f^{(n)}$):\\
By Lemma~\ref{lemfuestimates} (i) and (iii),
\beq \label{lasteq2}
\|\pa_{z}f^{(n)}(t)\|_\infty
\leq 
C\left(1+\int_0^t\big(1+
\|\pa^2_{(t,x)}u^{(n-1)}(s)\|_\infty\|\pa_{z}f^{(n)}(s)\|_\infty\big)ds\right),
\eeq
\beq \label{lasteq}
\|\pa^2_{(t,x)}u^{(n)}(t)\|_\infty
\leq
C\left((1+\log^*(\sup_{0\leq t\leq T_0}\|\pa_x f^{(n)}(t)\|_\infty)
+\int_0^t\|\pa^2_{(t,x)}u^{(n-1)}(s)\|_\infty ds\right).
\eeq
Applying Gronwall's lemma to (\ref{lasteq2}) yields
\beq\label{pafn}
 \|\pa_{z}f^{(n)}(t)\|_\infty
\leq 
C\exp\left(C\int_0^t(1+\|\pa^2_{(t,x)}u^{(n-1)}(s)\|_\infty)ds\right).
\eeq
Inserting this into (\ref{lasteq}) we obtain
\[
\|\pa^2_{(t,x)}u^{(n)}(t)\|_\infty
\leq 
C+C\int_0^t(1+\|\pa^2_{(t,x)}u^{(n-1)}(s)\|_\infty)ds,
\]
and by induction,
\[ 
 \|\pa^2_{(t,x)}u^{(n)}(t)\|_\infty\leq Ce^{CT_0}.
\] 
The bound on $\|\pa_{z}f^{(n)}(t)\|_\infty$ now follows 
from (\ref{pafn}).

\noindent
{\em Step 3 } (Uniform Cauchy property of $f^{(n)}(t)$ and $\pa_x u^{(n)}(t)$): \\
Introduce
\[
f^{m,n}= f^{(m)} - f^{(n)}, \qquad
u^{m,n}= u^{(m)} - u^{(n)}
\]
and note that
\[
S f^{m,n} = (\pa_x u^{(n-1)}) \pa_v f^{m,n} + \pa_x u^{m-1,n-1} \pa_v f^{(m)}.
\]
As in the uniqueness part we derive the estimates
\[
\| \pa_x u^{m,n}(t) \|_\infty \le C\int_0^t (\| \pa_x u^{m-1,n-1}(s) \|_\infty
+ \| f^{m,n}(s) \|_\infty) ds
\]
and
\[
\| f^{m,n}(t) \|_\infty \le C \int_0^t \| \pa_x u^{m-1,n-1}(s) \|_\infty ds.
\]
Combining these we obtain
\[
\| \pa_x u^{m,n}(t) \|_\infty \le C \int_0^t \| \pa_x u^{m-1,n-1}(s) \|_\infty ds,
\]
and hence by induction,
\[
\| \pa_x u^{m,n}(t) \|_\infty \le C \frac{(C t)^k}{k!}, \quad k=\min(m,n).
\]
So $\pa_x u^{(n)}(t)$ is a uniform Cauchy sequence, 
and the same is true for $f^{(n)}(t)$.

\noindent{\em Step 4} (Uniform Cauchy property of 
$\pa_{z} f^{(n)}(t)$ and $\pa_{(t,x)}^2 u^{(n)}(t)$): \\
We begin by establishing a bound on 
$\|\pa_{(t,x)}^2 u^{m,n}(t)\|_\infty$. Using the representation of 
Lemma~\ref{lem2derivative}, proceeding on a term-by-term basis and again using
a splitting as in the uniqueness part we obtain
\beq\label{est2u}
 \|\pa_{(t,x)}^2 u^{m,n}(t)\|_\infty\leq 
 C \int_0^t\left(\|\pa_{(t,x)}^2 u^{m-1,n-1}(s)\|_\infty
            +\|\pa_{z}f^{m,n}(s)\|_\infty \right)\, ds.
\eeq
Here the critical $TT$-kernel is estimated as in the proof of 
Lemma~\ref{lemfuestimates} letting $\tau\to t$. For the other terms
we use the boundedness assertions as well as the Cauchy properties
already obtained for the lower derivatives. 

Next we prove that the characteristics converge uniformly. 
Writing $Z^{(n)}(s,t,z)=(X^{(n)}(s,t,z),V^{(n)}(s,t,z)),\ 
z=(x,v)$, where we omit the arguments if there is no danger of 
misinterpretation, we find
\beas
\left|\frac{d}{ds}(X^{(n)}-X^{(m)})\right|
&\leq&
|V^{(n)}-V^{(m)}|, \\  
\left|\frac{d}{ds}(V^{(n)}-V^{(m)})\right|
&=&
|\pa_x u^{(n-1)}(X^{(n)})-\pa_x u^{(m-1)}(X^{(m)})|\\
&\leq&
|\pa_x u^{(n-1)}(X^{(n)})-\pa_x u^{(n-1)}(X^{(m)})|\\
&&
\qquad +|\pa_x u^{(n-1)}(X^{(m)})-\pa_x u^{(m-1)}(X^{(m)})|\\
&\leq&
C\,|X^{(n)}-X^{(m)}|+\delta_{m,n}.
\eeas
Here the expression $\delta_{m,n}$ converges to zero if $m,n\to\infty$ by the
Cauchy property of $\pa_x u^{(n)}$ and we have used the boundedness
of $\pa^2_xu^{(n)}$. Combining these two estimates and again using Gronwall's
lemma we obtain the claimed convergence of the characteristics, i.e.,
$Z^{(n)}(s,t,z)$ converges uniformly for $0\leq s\leq T_0$; the convergence
also uniform w.r.t.\ the parameters $t$ and $z$.

Writing $Z^{(n)}(s)$ for $Z^{(n)}(s,t,z)$ the analog of equation (\ref{inteqfz}) 
for iterates and vanishing $g$ implies
\beas
\pa_z f(t,z) 
&=&
\pa_z\fn (Z^{(n)}(0))-
    \int_0^t \pa_z f^{(n)}(s,Z^{(n)}(s))\ \pa^2_xu^{(n-1)}(s,X^{(n)}(s))\,ds,
\eeas
and hence
\beas
 |\pa_z f^{m,n}(t,z)|
 &\leq& |\pa_z\fn(Z^{(m)}(0))-\pa_z\fn(Z^{(n)}(0))|\\
 &&+\int_0^t\Bigl|\pa_z f^{(m)}(s,Z^{(m)}(s))\ \pa^2_xu^{(m-1)}(s,X^{(m)}(s))\\
 &&\hspace*{2.5cm}
   -\pa_z f^{(n)}(s,Z^{(n)}(s))\ \pa^2_xu^{(n-1)}(s,X^{(n)}(s))\Bigr|\, ds.
\eeas
The first term vanishes in the limit $m,n\to\infty$, and we split the second term
into the following four parts:
\bea\label{problem}\nn
&&\int_0^t\Big(
    |\pa_z f^{(m)}(s,Z^{(m)}(s))-\pa_z f^{(m)}(s,Z^{(n)}(s))|\ 
    |\pa^2_xu^{(m-1)}(s,X^{(m)}(s))|\\
  &&\qquad+|\pa_z f^{(m)}(s,Z^{(n)}(s))|\ 
    |\pa^2_{x}u^{(m-1)}(s,X^{(m)}(s))-\pa^2_{x}u^{(m-1)}(s,X^{(n)}(s))| \nn \\
  &&\qquad+|\pa_z f^{(m)}(s,Z^{(n)}(s))-\pa_z f^{(n)}(s,Z^{(n)}(s))|\ 
    |\pa^2_xu^{(m-1)}(s,X^{(n)}(s))|\nn\\
  &&\qquad+|\pa_z f^{(n)}(s,Z^{(n)}(s))|\ 
   |\pa^2_{x}u^{(m-1)}(s,X^{(n)}(s))-\pa^2_{x}u^{(n-1)}(s,X^{(n)}(s))|
  \Big)\ ds.\ \ \ \ \ \ 
\eea
The first two vanish in the limit $m,n\to\infty$ by Lemma~\ref{missing} below
while the latter two may be estimated due to the boundedness of $\pa^2_xu^{(n)}$
and $\pa_z f^{(n)}$. Hence
\[
\|\pa_z f^{m,n}(t,z)\|_\infty
\leq
\delta_{m,n}+C\int_0^t\left(\|\pa_z f^{m,n}(s)\|_\infty
+\|\pa^2_xu^{m-1,n-1}(s)\|_\infty \right)\,ds,
\]
where again $\delta_{m,n}\to0$ if $m,n\to\infty$. Another application
of Gronwall's lemma gives
\beq\label{fcauchy}
  \|\pa_z f^{n,m}(t,z)\|_\infty
  \leq\delta_{m,n}+C\int_0^t\|\pa^2_xu^{m-1,n-1}(s)\|_\infty ds.
\eeq
We finally insert the last inequality into equation (\ref{est2u})
to obtain
\[ 
 \|\pa_{(t,x)}^2u^{m,n}\|_\infty
  \leq\delta_{m,n}+C\int_0^t\|\pa^2_{(t,x)}u^{m-1,n-1}(s)\|_\infty ds
\]
and, by iteration for each $l\in\N$
\[
  \|\pa_{(t,x)}^2u^{m,n}\|\infty\leq\delta_{m,n}\ e^{C}+\frac{BC^lT_0^l}{l!},
\]
on $[0,T_0]$, where $B$ is a bound for $\|\pa_{(t,x)}^2u^{m}\|_\infty$. 
This proves the required
Cauchy property for $\|\pa_{(t,x)}^2u^{n}\|_\infty$ while the one for 
$\|\pa_z f^n(s)\|_\infty$ follows from equation (\ref{fcauchy}),
thereby finishing the proof.
\end{proof}

It remains to show that the first two terms in (\ref{problem})
indeed vanish if $m,n\to\infty$, which amounts to showing that $\pa_z f^{(n)}$ as
well as $\pa^2_{(t,x)}u^{(n)}$ are uniformly continuous, with modulus of continuity
uniform in $n$. One
could be tempted to use the mean value theorem but suitable 
estimates on the second order derivatives of $f^{(n)}$ and the third 
order derivatives of $u^{(n)}$ are not available. We begin by defining 
\beas
\delta_n(t,\eta)
&:=&
\sup\left\{|\pa_z f^{(n)}(t,z)-\pa_z f^{(n)}(t,z')|
\mid z,z'\in\R^6,\ |z-z'|\leq\eta \right\}, \\
\theta_n(t,\eta)
&:=&
\sup\left\{|\pa^2_{(t,x)}u^{(n-1)}(t,x)-\pa^2_{(t,x)}u^{(n-1)}(t,y)|
\mid x,y\in\R^3,\ |x-y|\leq\eta\right\}.    
\eeas
Note that both $\delta_n$ and $\theta_n$ are bounded uniformly in $n$.
Our desired result now is

\begin{lemma}\label{missing}
On any time interval $[0,T_0]$ on which the iterates satisfy the bounds
established in Steps~1 and 2 in the proof of Thm.~\ref{locsol} the
following is true:
For all $\varepsilon>0$ there exist $\eta_0>0$ and $n_0\in \N$ such that for all 
$n\geq n_0$
\[
\delta_n(t,\eta_0),\ \theta_n(t,\eta_0)\leq\varepsilon,\ t\in[0,T_0]. 
\]
\end{lemma}
We defer the rather technical proof of this lemma to Appendix~\ref{ul}.
We now establish a continuation criterion for the solutions
obtained in Thm.~\ref{locsol}:

\begin{theorem}\label{contcrit}(Continuation criterion)
 Let $(f,u)$ be a solution of the relativistic Vlasov-Klein-Gordon system
 on $[0,T[$ as in Theorem~\ref{locsol}. Then the function
\[
P(t):=\sup\left\{|v|\ | f(\tau,x,v)\not=0,\ 0\leq\tau\leq t,\ x\in\R^3\right\}
\]
is bounded on $[0,T[$ iff $\|\pa_{x}u(t)\|_\infty$
is bounded on $[0,T[$.
Moreover, if $T$ is chosen maximally then any of these bounds
implies that the solution is global, i.e., $T=\infty$.
\end{theorem}

\begin{proof}
To prove that a bound on $\pa_{x}u$ implies a bound on $P$
we integrate the $v$-component of the characteristic
system.
For the reverse direction we note that we can estimate $\pa_{x}u$
in terms of $P$ exactly as we did for the iterates in Step~1
of the proof of Thm.~\ref{locsol}, using 
Lemma~\ref{lem1derivative}. 

Assume now that $T$ is chosen maximally, that $P$ is bounded on
$[0,T[$, and $T<\infty$. For any $t_0 \in [0,T[$ we can use the arguments from the
proof of Thm.~\ref{locsol} to show that a solution with data
$f(t_0),\ u(t_0),\ \pa_t u(t_0)$ prescribed at $t=t_0$ exists
on some time interval $[t_0,t_0 + \delta[$, except that there is
one technical catch here: $u(t_0),\ \pa_t u(t_0)$
are not sufficiently regular to qualify as initial data in the context of
Thm.~\ref{locsol}. But since we already have the solution
on $[0,t_0]$ we can define the iterates used to obtain the extended solution
as follows: For $(f^{0}, u^{0})$ we take a global extension of the existing
solution with the required regularity and with $||f^{(0)}(t)||_\infty$,
$||\pa_z f^{(0)}(t)||_\infty$,
$||\pa_{(t,x)} u^{(0)}||_\infty$, $||\pa^2_{(t,x)} u^{(0)}||_\infty$
bounded in $t$. Given the $(n-1)$st iterate we define the $n$th iterate
exactly as before on $[0,\infty[$. Then all these iterates coincide
with the solution on $[0,t_0]$, the data term in the formulas for the field,
from which the loss of derivatives arises, is the one determined by
the data at $t=0$, and it is straight forward to repeat
the arguments from the proof of Thm.~\ref{locsol} to extend the solution to
some time interval $[0,t_0 + \delta[$.
The crucial point now is that the uniform bound
on the momenta implies that $\delta >0$ can be chosen independently of $t_0$,
cf.\ (\ref{pngronwall}) and the lines that follow.
For $t_0$ close enough to $T$ this contradicts the maximality of $T$.
\end{proof}

\section{The one-dimensional case} \label{1d}

In this section we illustrate that our continuation criterion
from Thm.~\ref{contcrit} holds and hence classical solutions are global
in the one dimensional case, where $x, v \in \R$. To do so we first need to
derive the representation formulas for $u$ and its derivatives in this situation.
The standard trick to do this is to observe that $u(t,x)$ solves the Klein-Gordon 
equation (\ref{kg}) iff $w(t,x, \xi) = u(t,x) e^{-i\xi}$ solves the wave equation
\[
(\pa^2_t - \pa^2_x  - \pa^2_\xi) w = -e^{-i\xi} \rho 
\]
with initial data transformed accordingly. This leads to
\[
u(t,x) = u_{\mathrm{hom}}(t,x) + u_{\mathrm{inh}}(t,x),\quad t\geq 0,\ x\in \R,
\]
where $u_{\mathrm{hom}}(t,x)$ is the solution of the homogeneous equation
and depends only on the initial data for $u$, and
\[
u_{\mathrm{inh}}(t,x) 
= -\frac{1}{2} \int_0^t \int_{x-(t-s)}^{x+(t-s)} \rho(s,y)\, 
J_0(\sqrt{(t-s)^2 - |x-y|^2})\, dy\, ds.
\]
Hence
\bea \label{1dfield}
\pa_x u_{\mathrm{inh}}(t,x) 
&=&
-\frac{1}{2} \int_0^t\left(\rho(s,x+(t-s)) - \rho(s,x-(t-s))\right)\, ds \nonumber \\
&&
{} - \frac{1}{2} \int_0^t \int_{x-(t-s)}^{x+(t-s)} \rho(s,y)\,
\frac{J_1(\sqrt{(t-s)^2 - |x-y|^2})}{\sqrt{(t-s)^2 - |x-y|^2}} (x-y)\, dy\, ds.\ \ 
\eea
Assume now that we have a (local) solution of the Vlasov-Klein-Gordon system
in the one-dimensional case, with $f(t)$ compactly supported for all $t$.
As before, 
\[
P (t):=\sup\left\{|v| \mid
f (\tau,x,v)\not=0,\ 0\leq\tau\leq t,\ x\in\R\right\}.
\]
Since $f$ is constant along characteristics of the Vlasov equation,
\[
\rho (t,x) \leq 2 ||\fn||_\infty P(t).
\]
Integrating the Vlasov equation w.r.t.\ $x$ and $v$ implies that
$||f(t)||_1 = ||\rho(t)||_1 = ||\fn||_1$,
and since $J_1(\xi)/\xi$ is bounded on $\xi >0$, (\ref{1dfield}) implies that
\[
||\pa_x u(t)||_\infty \leq C ( 1+ t^2 + t P(t)).
\]
Integrating the $v$-component of the characteristic system implies that
\beq \label{1dpgronwall}
P (t) \leq \Pn + C (1+t)^3 + C\int_0^t s P(s)\, ds
\eeq
so that by Gronwall's lemma $P$ can not blow up on any bounded time interval.

Indeed, these estimates could be repeated for iterates defined as in the
proof of Thm.~\ref{locsol}. Controlling first order derivatives of $f^{(n)}$
and second order ones of $u^{(n)}$ would be much easier than in the three dimensional
case, as should be obvious from comparing Eqn.~(\ref{1dfield}) with 
Lemma~\ref{lem1derivative}. Hence:

\begin{theorem}\label{1dglobal}
Let $\fn\; \in \CC_c^1(\R^2)$, $\uno \in \CC_b^3(\R)$,
$\unt \in \CC_b^2(\R)$.
Then there exists a unique classical solution 
\[ 
f \in \CC^1([0,\infty[ \times \R^2),\ u\in \CC^2([0,\infty[\times\R)
\]
of the one-dimensional relativistic Vlasov-Klein-Gordon system,
satisfying the initial conditions (\ref{id}). If $P$ denotes the solution of the
integral equation corresponding to (\ref{1dpgronwall}) then
\[
f(t,x,v) = 0 \ \mbox{for}\ |x| \geq \Rn + t\ \mbox{or}\ |v| \geq P(t)
\]
with $\Rn$ determined by $\fn$.
\end{theorem}

\begin{appendix}
\section{Some facts on the Vlasov equation}

In this appendix we collect some useful facts on the (inhomogeneous) Vlasov
equation for easy reference.
As before, we combine $x$ and $v$ to one variable $z=(x,v)$, and we consider
the initial value problem
\[
\pa_t f(t,z) + G(t,z) \pa_z f(t,z) = g(t,z), \qquad f(0,z)= \fn(z),
\]
where $\fn\;, G\; g\in \CC^1$ and $G$ is such that the solutions of the
corresponding characteristic system
\[
\dot{z}(s) = G(s,z(s))
\]
exist on the time interval on which $G$ is defined. 
Denote by $s \mapsto Z(s,t,z)$ the solution corresponding to the
initial condition $Z(t,t,z)=z$ and recall that
\[
Z(s,t,z)= Z(s,r,Z(r,t,z))
\]
and
\[
\pa_t Z(s,t,z) + \pa_z Z(s,t,z) G(t,z) = 0;
\]
the second equation follows from the first by differentiating $z_0= Z(s,t,Z(t,s,z_0))$
with respect to $t$ and then choosing $z_0=Z(s,t,z)$. The partial derivative
$\pa_z Z(s,t,z)$ satisfies the first variational equation
\[
\pa_z\dot{Z}(s,t,z) = (\pa_z G)(s, Z(s,t,z)) \, \pa_z Z(s,t,z), 
\]
with $\pa_z Z(s,s,z) = \mathbb{I}$ the unit matrix, or equivalently,
\[
\pa_z Z(s,t,z) = \mathbb{I} + \int_t^s (\pa_z G)(r, Z(r,t,z)) \pa_z Z(r,t,z) dr.
\]
In addition, there is also the less obvious equation
\begin{equation}\label{rol}
\pa_z Z(s,t,z) = \mathbb{I} + \int_t^s  \pa_z Z (s,r,Z(r,s,z)) (\pa_z G)(r, Z(r,t,z))\,dr,
\end{equation}
which holds because the right hand side solves
the first variational equation. However, note that the integrands in the last two
equations are not equal.

We now apply these results to the Vlasov equation. Clearly the solution
is given by
\begin{equation} \label{solvl}
f(t,z)= \fn(Z(0,t,z)) + \int_0^t g(s,Z(s,t,z))\,ds.
\end{equation}
Moreover, $\pa_z f(t,z)$ exists and satisfies
\bea  \nn
\pa_z f (t,z) &=& (\pa_z\fn) (Z(0,t,z)) - \int_0^t (\pa_z f)(s,Z(s,t,z) 
(\pa_z G)(s,Z(s,t,z)) ds\\ \label{inteqfz}
&& {} + \int_0^t (\pa_z g)(s,Z(s,t,z)\,ds.
\eea
This can be seen by a straightforward calculation using equation~(\ref{rol}).
We should remark that the usual way of deriving this equation
is by differentiating the Vlasov equation with respect to $z$ and viewing the resulting
equation as an inhomogeneous Vlasov equation for $\pa_z f$. This
requires $\fn\:, g \in \CC^2$ which was not necessary here.

\section{Proof of Lemma~\ref{missing}}\label{ul}

We split the proof into several steps, following the approach used
in \cite[pp.\ 78--90]{K} in the case of the Vlasov-Maxwell system.
Throughout we argue on a time interval $[0,T_0]$ on which the
iterates satisfy all the bounds established in Steps~1 and 2 of 
the proof of Thm.~\ref{locsol}.

\begin{lemma}\label{lemA}
(Estimate on $\theta_n$ in terms of $\theta_{n-1}$, $\delta_{n-1}$)\\
There exists $C>0$ such that $\forall\varepsilon>0\ \exists\eta_0>0\
\forall \eta\in[0,\eta_0]\ \forall t\in[0,T_0]\ \forall n\in\N$,
\[
\theta_n(t,\eta)
\leq
\eps + C\, \left(\eta+\int_0^t\big(\theta_{n-1}(s,\eta)
+\delta_{n-1}(s,\eta)\big)\,ds \right).
\]
\end{lemma}

\begin{proof}
Using once more the representation formulas of Lemma~\ref{lem2derivative}
we proceed as in Steps~1 and 2 of the proof of Theorem~\ref{locsol}. 
\end{proof}

\begin{lemma}\label{lemF}
(Estimate on $\delta_n$ in terms of $\theta_{n}$)\\
There exists $C >0$ such that $\forall\varepsilon>0\ 
 \exists\eta_0>0\
 \forall \eta\in[0,\eta_0]\ \forall t\in[0,T_0]\ \forall n\in\N$,
 \[
  \delta_n(t,\eta)\leq\eps+C\int_0^t\theta_{n}(s,\eta)ds.
 \]
\end{lemma}

To prove Lemma~\ref{lemF} we need some additional
information on the derivatives of the characteristics which is provided by
the following lemma.

\begin{lemma}\label{subchar} 
(Estimates on the derivatives of the characteristics)
\begin{itemize}
\item [(i)]  
There exists $C>0$ such that $\forall s, t \in[0,T_0]$ 
$\forall z,z'\in\R^6\ \forall n\in\N$
\[
|Z^{(n)}(s,t,z)-Z^{(n)}(s,t,z')|\leq C|z-z'|
\]
\item [(ii)] 
There exists $C>0$ such that $\forall s, t\in[0,T_0]$, $s\leq t$,
$\forall z,z'\in\R^6\ \forall n\in\N$
\beas
&&
\quad|\pa_z Z^{(n)}(s,t,z)-\pa_z Z^{(n)}(s,t,z')|\\
&&
\qquad
\leq C|z-z'|+\int_s^t\left|\pa^2_xu^{(n-1)}(\tau,X^{(n)}(\tau,t,z))-
                        \pa^2_xu^{(n-1)}(\tau,X^{(n)}(\tau,t,z'))\right|\,d\tau.
\eeas
\end{itemize}
\end{lemma}
\begin{proof}
Part (i) follows immediately from equation (\ref{rol}).
As to (ii), we only treat the terms involving $x$-derivatives; the
$\pa_v$-terms may be estimated in the same fashion. So we start
with the term $|\pa_{x_i}X^{(n)}(s,t,z)-\pa_{x_i}X^{(n)}(s,t,z')|$, 
$1\leq i \leq 3$.
By an elementary calculation, cf.\ \cite[pp.\ 83--85]{K}, we obtain
\bea\label{*ECT5}
&&
|\pa_s\pa_{x_i}X^{(n)}(s,t,z)-\pa_s\pa_{x_i}X^{(n)}(s,t,z')| \nn \\
&&
\qquad \leq
2|\pa_{x_i}V^{(n)}(s,t,z)-\pa_{x_i}V^{(n)}(s,t,z')|\nn\\
&&\qquad \quad {}+5|\pa_{x_i}V^{(n)}(s,t,z)-\pa_{x_i}V^{(n)}(s,t,z')|\
|V^{(n)}(s,t,z)-V^{(n)}(s,t,z')|.\ \ \ 
\eea
Using 
\[
\pa_s\pa_{x_i}V^{(n)}(s,t,z)=-\pa_{x_i}\pa_x u^{(n-1)}(s,X^{(n)}(s,t,z))\ 
\pa_{x_i}X^{(n)}(s,t,z)
\]
we obtain
\bea\label{**ETC6}  
&&
|\pa_s\pa_{x_i}V^{(n)}(s,t,z)-\pa_s\pa_{x_i}V^{(n)}(s,t,z')| \nn \\
&&
\qquad\leq
|\pa^2_xu^{(n-1)}(s,X^{(n)}(s,t,z))-\pa^2_xu^{(n-1)}(s,X^{(n)}(s,t,z'))|\ 
|\pa_{x_i}X^{(n)}(s,t,z)| \nn \\
&&
\qquad \quad {}+|\pa^2_xu^{(n-1)}(s,X^{(n)}(s,t,z'))|\ 
|\pa_{x_i}X^{(n)}(s,t,z)-\pa_{x_i}X^{(n)}(s,t,z')|.
\eea
Using (i) and equations (\ref{*ECT5}) and (\ref{**ETC6}) 
we obtain
\beas
&&
|\pa_{x_i}Z^{(n)}(s,t,z)-\pa_{x_i}Z^{(n)}(s,t,z')|\\
&&
\qquad \leq C\Big(|z-z'|
+\int_s^t |\pa_{x_i}Z^{(n)}(\tau,t,z)-\pa_{x_i}Z^{(n)}(\tau,t,z')|d\tau\\
&&
\qquad \hphantom{\leq C\big(|z-z'|}  
+\int_s^t|\pa^2_xu^{(n-1)}(\tau,X^{(n)}(\tau,t,z))-\pa^2_xu^{(n-1)}   
(\tau,X^{(n)}(\tau,t,z'))|d\tau\ \Big).
\eeas
Another appeal to Gronwall's lemma completes the proof of the lemma.
\end{proof}

\begin{proof}[Proof of Lemma~\ref{lemF}]
  We only treat the $\pa_v f1{(n)}$-terms; the $\pa_x f^{(n)}$-terms can be estimated 
  analogously. Again writing $z=(x,v)$ respectively $z'=(y,w)$ we find 
  \bea\label{I-IV}\nn
   &&\hspace*{-1cm}|\pa_v f^{(n)}(t,z)-\pa_v f^{(n)}(t,z')|\\
   &\leq&
   |\pa_z\fn(Z^{(n)}(0,t,z))|\ |\pa_vZ^{(n)}(0,t,z)-\pa_vZ^{(n)}(0,t,z')| \nn \\
   &&+|\pa_v Z^{(n)}(0,t,z')|\ 
     |\pa_z\fn(Z^{(n)}(0,t,z))-\pa_z\fn(Z^{(n)}(0,t,z'))|.
  \eea
  The second term in the above estimate is bounded by $\eps$ for
  large $n$, $t\in [0,T_0]$ and $|z-z'|$ suitably small by the fact 
  that $\pa_z Z^{(n)}$ is uniformly (in $n$) bounded, 
  the uniform continuity of the data $\pa_z\fn$ and Lemma
  \ref{subchar} (i).

  By Lemma \ref{subchar} (ii) and the fact that the characteristics are
  uniformly bounded the first term in equation (\ref{I-IV})
  may be estimated by
  \[
   C_1 |z-z'| + C_1\int_0^t\left|\pa^2_xu^{(n-1)}(\tau,X^{(n)}(\tau,t,z))-\pa^2_xu^{(n-1)}
   (\tau,X^{(n)}(\tau,t,z'))\right|\,d\tau.
  \]
  Denote by $C_2$ the maximum of $C_1$ and $C$ from Lemma \ref{subchar} (i).
  Set $C=(\lceil 2C_2\rceil+1)C_1$ where
$\lceil r \rceil$ denotes the
smallest integer bigger or equal to $r$. We claim that this constant verifies
  the assertion. Indeed, let $\eps>0$, and $\eta\in [0,\eta_0]$ with 
  $\eta_0:=\eps/(2C_1+1)$. Let $z,z'\in\R^6$ with $|z-z'|\leq\eta$, $n\in\N$ and
  $\tau,t\in[0,T_0]$. Then by Lemma \ref{subchar} (i) 
  \[
   |Z^{(n)}(\tau,t,z)-Z^{(n)}(\tau,t,z')|\leq 2C_2\eta.
  \]
  Together with the fact that $\theta_n(\tau,k\eta)=k\theta_n(\tau,\eta)$
  for all $k\in \N$ this gives
  \beas
    &&\hspace*{-2cm}
    |\pa^2_xu^{(n-1)}(\tau,X^{(n)}(\tau,t,z))-
             \pa^2_xu^{(n-1)}(\tau,X^{(n)}(\tau,t,z'))|\\
    &\leq&\theta_n(\tau,2C_2\eta) \leq \theta_n(\tau,(1+\lceil2C_2\rceil)\eta)
    \leq (1+\lceil2C_2\rceil)\theta_n(\tau,\eta).
  \eeas
  Summing up we obtain
\[
|\pa_v f^{(n)}(t,z)-\pa_v f^{(n)}(t,z')|
\leq
\eps+C \int_0^t\theta_n(\tau,\eta)\,d\tau.
\]
\end{proof}

\noindent
{\em Proof of Lemma~\ref{missing}.}
 Combining Lemma~\ref{lemA} and Lemma~\ref{lemF} we obtain: 
\[
\forall\eps>0\ 
\exists\eta_0>0\ \forall\eta\in[0,\eta_0],\ t\in[0,T_0],\ n\in \N:\;
\theta_n(t,\eta)\leq C\left(\eps+\int_0^t\theta_{n-1}(s,\eta)ds\ \right),
\]
which by iteration shows that 
\[
\forall \eps>0\ \exists\eta_0>0\
\exists n_0\in\N\ \forall n\geq n_0,\ t\in [0,T]:\
\theta_n(t,\eta_0)\leq\eps.
\]
The claim concerning $\delta_n$ now immediately follows from Lemma~\ref{lemF}.
\ep

\end{appendix}


\begin{thebibliography}{10}

\bibitem{Ab}
Abraham, M.,
{\em Theorie der Elektrizit\"at, Band 2: Elektromagnetische Theorie der
Strahlung}.
Teubner, Leipzig 1905.

\bibitem{DL}
DiPerna, R. J., Lions, P.-L.,
Global weak solutions of {V}lasov-{M}axwell systems.
{\em Comm.\ Pure Appl.\ Math.}\ {\bf 42}(6), 729--757 (1989).

\bibitem{GlStr2}
Glassey, R. T.,  Strauss, W. A.,
Singularity formation in a collisionless plasma could occur only at high velocities.
{\em Arch.\ Rat.\ Mech.\ Anal.}\ 
{\bf 92}, 59--90 (1986).

\bibitem{GlStr}
Glassey, R. T.,  Strauss, W. A.,
Absence of shocks in an initially dilute collisionless plasma.
{\em Commun.\ Math.\ Phys.}\ 
{\bf 113}, 191--208 (1987).

\bibitem{IKM}
Imaikin, V.~M., Komech, A.~I., Markowich, P.~A.,
Scattering of solitons of the Klein-Gordon equation coupled to a classical
particle. {\em J. Math. Phys.}\ {\bf 44}(3), 1202--1217 (2003).

\bibitem{IKS}
Imaikin, V.~M., Komech, A.~I., Spohn, H.,
Scattering theory for a particle coupled to a scalar field.
{\em Discrete Contin. Dyn. Syst.}\ {\bf 10}(1-2), 387--396 (2004).

\bibitem{KKS1}
Komech, A.~I., Kunze, M., Spohn, H.,
Long-time asymptotics for a classical particle coupled to a scalar wave field.
{\em Commun.\ Part.\ Diff.\ Equations} {\bf 22}, 307--335 (1997).

\bibitem{KKS2}
Komech, A.~I., Kunze, M., Spohn, H.,
Effective dynamics for a mechanical particle coupled to a wave field.
{\em Commun.\ Math.\ Phys.} {\bf 203}, 1--19 (1999).

\bibitem{KS}
Komech, A.~I., Spohn, H.,
Soliton-like asymptotics for a classical particle interacting
with a scalar wave field.
{\em Nonlin.\ Anal.} {\bf 33}, 13--24 (1998).

\bibitem{K}
Kunze, M., Das relativistische Vlasov-Maxwell System partieller
Differentialgleichungen zu Anfangsdaten mit nichtkompaktem Tr\"ager.
Diploma-thesis, University of Munich (1991).

\bibitem{KRST}
Kunzinger, M., Rein, G., Steinbauer, R., and Teschl, G.,
Global Weak Solutions of the Relativistic Vlasov-Klein-Gordon System.
{\em Commun.\ Math.\ Phys.}\ 
{\bf 238}, 367--378 (2003).

\bibitem{MS}
Morawetz, C.~S., Strauss, W.~A.,
Decay and scattering of solutions of a nonlinear relativistic wave
equation.
{\em Commun.\ Pure Applied Math.}\ {\bf 25}, 1--31 (1972).

\bibitem{R}
Rein, G.,
Global weak solutions to the relativistic Vlasov-Maxwell system
revisited. {\em Comm.\ Math.\ Sci.}\ {\bf 2 (2)}, 145--158 (2004),
in press.

\bibitem{Sid}
Sideris, T.~C.,
Decay estimates for the three-dimensional inhomogeneous
Klein-Gordon equation and applications.
{\em Commun.\ Part.\ Diff.\ Eqns.}\ {\bf 14}, 1421--1455 (1989).
 
\end{thebibliography}
\end{document}